\documentclass[a4paper,10.5pt]{article}
\usepackage{theorem}
\usepackage{latexsym,amssymb,amsfonts,amsmath}
\usepackage[dvips]{graphicx}
\usepackage{epsfig}
\usepackage{caption}

\newtheorem{definition}{Definition}
\newtheorem{theorem}{Theorem}
\newtheorem{lemma}{Lemma}

\newtheorem{proposition}[theorem]{Proposition}

\usepackage{latexsym}

\title{{\Large {\bf On a property of the simple random walk on $\mathbb{Z}$ 
}
}}

\author{ 
{\small 
Norio Konno,$^{1}$ 
\footnote{konno@ynu.ac.jp}
\quad Hayato Saigo$^{2}$
\footnote{h\_saigoh@nagahama-i-bio.ac.jp}
and Hiroki Sako $^{3}$
\footnote{sako@ie.niigata-u.ac.jp}
\quad  
}\\ 
{\scriptsize $^{1}$ 
Department of Applied Mathematics, Faculty of Engineering, Yokohama National University,,
}\\
{\scriptsize Hodogaya, Yokohama, 240-8501, Japan
} \\
{\scriptsize $^2$ 
Nagahama Institute of Bio-Sciences and Technology
}\\
{\scriptsize Tamura, Nagahama, 526-0829, Japan
} \\
{\scriptsize $^3$ 
Faculty of Engineering, Niigata University,
}\\
{\scriptsize Nishi-ku, Niigata 950-2181, Japan
} \\
} 

\date{\empty }
\begin{document}
\maketitle

\par\noindent
\begin{small}
\baselineskip=24pt
\par\noindent
{\bf Abstract}
The subject of this paper is the simple random walk on $\mathbb{Z}$.
We give a very simple answer to the following problem: under the condition 
that a random walk has already spent $\alpha$-percent 
of the traveling time on the positive side $\mathbb{Z}_{\ge 0}$, 
what is the probability that the random walk is now on the positive side?

The symmetric random walks which step $2n$-times can be decomposed 
in the following two ways: (1) how many times the walk steps on the positive side, 
(2) whether the last step is on the positive side or on the negative side.
To answer the problem above, we clarify the number of the walks classified by (1) and (2).
It has been already known that the distribution of the number indicated by (1) makes the arcsine law.
Combining with the decomposition with respect to (2), we obtain a decomposition of the arcsine law into the Marchenko-Pastur law.

\footnote[0]{{\it Key words. Random walk, Arcsine law, Marchenko-Pastur law} }

\end{small}

\baselineskip=24pt

\setcounter{equation}{0}
\section{Introduction}

An invisible simple random walker on $\mathbb{Z}$ suddenly telephoned to you at the time $t=2n$ and informed you that she has spent $70$ persent of time on 
the positive side, from $t=0$ to $t=2n$.

Question: What is the probability that she is on the positive side \textbf{now} ($t=2n$)?

Answer: $70$ persent.

\dots But WHY? 

In the present paper, we will explain WHY (and to show some remarkable corollalry connecting classical probability theory and noncommutative probability). 
In Section 2 we briefly review basic definitons and known results on the simple random walk on $\mathbb{Z}$. 
We prove the main theorem (Theorem \ref{main theorem}) and as a corollary we explain WHY (Theorem \ref{why}) in Section 3.
In the last section we show that the theorem implies the new limit theorem connecting simple random walk with positive ends and the Marchenko-Pastur law.

\section{Preliminaries on the simple random walk on $\mathbb{Z}$}

\textbf{notation: discrete time interval $[m,n]:=\{m,m+1,\dots,n\}$ for integers satisfying $m\leq n$.}




\begin{definition}[$\mathcal{P}_N$,$\mathcal{P}_N^0$]

We denote the set of all possible paths of simple random walk $X(t)$ on discrete time interval $[0,N]=\{0,1,2,3,\dots ,N\}$ by $\mathcal{P}_N$.
The set of all possible paths of simple random walk $X(t)$ on time interval $[0,N]$ such that $X(N)=0$ is denoted by $\mathcal{P}_N^0$.

\end{definition}

These sets can be written as
\begin{eqnarray*}
\mathcal{P}_N
&=&
\{X \colon [0, N] \rightarrow \mathbb{Z} \ | 
\ X(0) = 0, |X(t) - X(t - 1)| = 1\ (t \in [1, N])\}, \\
\mathcal{P}_N^0
&=&
\{X \in \mathcal{P}_N \ | \ X(N) = 0\}.
\end{eqnarray*}
The number of the elements of $\mathcal{P}_N$ is $2^{N}$.
The number of the elements of $\mathcal{P}_N^0$ can be written as follows:

\begin{theorem}[Feller \cite{FEL}, Chap III Sec 4 Theorem1-(a)]\label{total}

\[
\#\mathcal{P}_N^0=\left\{
\begin{array}{ll}
 \binom{2n}{n}, &\quad \text{$N=2n$, \quad $(n=0,1,2,3,\dots)$}\\
 0, &\quad N {\rm\ is\ odd.}
\end{array}
\right.
\]

\end{theorem}

\begin{definition}[positive side]
Let $X(t)$ be a simple random walk on $\mathbb{Z}$. For positive integer $t$, 
we say that ``$X(t)$ is on the positive side'' if $X(t)\geq  0$ and $X(t-1)\geq 0$. In such a case, we also say that the path $X$ is on the positive side at $t$.
\end{definition}

\begin{definition}[$\mathcal{P}_N^{+}$]
We denote by $\mathcal{P}_N^{+}$ the set of all the random walks $X(t)$ defined for $t \in [0,N]$ which are on the positive side at $N$.
\end{definition}

\begin{definition}[soujourn time]
Let $X(t)$ be the simple random walk on $\mathbb{Z}$. The soujourn time on the positive side $T_{N}$ is defined by
\[
T_{N}:=\#\: \{t\in [1,N] \ |\ X(t) {\rm \ is\ on\ the\ positive\ side.}\}.
\]

\end{definition}

\begin{definition}[$\mathcal{P}_{N,m},\mathcal{P}_{N,m}^{0},\mathcal{P}_{N,m}^{+}.$]
We define $\mathcal{P}_{N, m}$ by
\[
\mathcal{P}_{N,m}=\{p\in \mathcal{P}_N \ |\ 
p {\rm \ is\ a\ path\ for\ which\ } T_{N}=m {\rm\ is\ satisfied. } \}.
\]
We also define $\mathcal{P}_{N,m}^{0}:=\mathcal{P}_{N,m}\cap \mathcal{P}_N^0 $ and $\mathcal{P}_{N,m}^{+}:=\mathcal{P}_{N,m}\cap \mathcal{P}_N^{+} $.
\end{definition}
For $\#\mathcal{P}_{2n,2k}$ and $\#\mathcal{P}_{2n,2k}^{0}$, the following general formula is well-known.

\begin{theorem}[Arcsine law for simple random walk, Feller \cite{FEL}, Chap III Sec 5 Theorem 1]\label{arcsine}
\[
\#\mathcal{P}_{2n,2k}=\binom{2k}{k}\binom{2n-2k}{n-k}.
\]
\end{theorem}

\begin{theorem}[Uniform principle, Feller \cite{FEL}, Chap III Sec 2 Theorem 3]\label{uniform}
\[
\#\mathcal{P}_{2n,2k}^{0}=\frac{1}{n+1}\binom{2n}{n}, 
\quad \mathrm{for\ every\ } k \in [0, n].
\]

\end{theorem}
Then how about $\#\mathcal{P}_{2n,2k}^{+}$? 
In the case that $k = n$, 
it is obvious that $\#\mathcal{P}_{2n,2n}^{+}=\#\mathcal{P}_{2n,2n}$, so Theorem \ref{arcsine} gives the following:

\begin{proposition}[Feller \cite{FEL}, Chap III Sec 4 Theorem 1-(c)]\label{positive}
\[
\#\mathcal{P}_{2n,2n}^{+}=\binom{2n}{n}.
\]
\end{proposition}
Theorem \ref{main theorem} in the next section gives a general formula for $\#\mathcal{P}_{2n,2k}^{+}$. As a collorary (Theorem \ref{why}), we get an answer to the question ``WHY'' in the introduction.

\section{Main Theorem}

\begin{lemma}\label{path counting}

\[
\#\mathcal{P}_{2n,2k}^{+}=\sum_{l\in[1,k]}\frac{1}{n-l+1}\cdot \binom{2n-2l}{n-l}\cdot 2\binom{2l-2}{l-1}.
\]

\end{lemma}
\textbf{Proof.}
Since the proposition is trivial for the case of $k=0,n$, we focus on the case of 
$k \in [1, n - 1]$.

It is obvious that for $X \in \mathcal{P}_{2n,2k}^{+}$
\[
\{t\in [1,2n-1] \ |\  X(t)=0 \}\neq \emptyset .
\]

Let $2\tau := {\rm Max}\{t\in [1,2n-1] \ |\  x(t)=0 \}$. Then $\tau \in [n-k,n-1]$. 
We can totally decompose the set $\mathcal{P}_{2n,2k}^{+}$ by the values of $2\tau$. It is easy to 
see that $X(t)$ is on the positive side for $2k-(2n-2\tau)$ times in $[0,2\tau]$ because  $X(t)$ is always (i.e. for $2n-2\tau$ times) 
on the positive side; and more strongly,

\begin{itemize}
 \item $X(2\tau +1)=1$
 \item $X(t)\geq 1$ on $[2\tau +1, 2n-1]$
 \item $X(2n)=X(2n-1)+1$ or $X(2n-1)-1$.
\end{itemize}

It is easy to see that the number of paths defined on $[2\tau+1,2n]$ satisfying 
$X(2\tau +1)=1$ and $X(t)\geq 1$ for all $t\in [2\tau+1,2n-1]$ is nothing but $2 \cdot \#\mathcal{P}_{2n-2\tau-2,2n-2\tau-2}^{+}$. Hence, from the above argument and Theorem \ref{uniform}, we obtain 

\[
\#\mathcal{P}_{2n,2k}^{+}=\sum _{\tau\in[n-k,n-1]}\#{P}_{2\tau,2k-(2n-2\tau)}^{0}\cdot (2 \cdot \#\mathcal{P}_{2n-2\tau-2,2n-2\tau-2}^{+})
\]
\[
=\sum _{\tau\in[n-k,n-1]}\frac{1}{\tau +1}\binom{2\tau}{\tau}\cdot 2\binom{2(n-\tau)-2}{(n-\tau)-1}
\]
\[
=\sum_{l\in [1,k]}\frac{1}{n-l+1}\binom{2n-2l}{n-l}\cdot 2\binom{2l-2}{l-1},
\]
where $l:=n-\tau$. $\Box$

\begin{lemma}\label{identity}
Let $k,l,n$ be positive integers and $1\leq k \leq n-1$. Then the identity
\[
\sum_{l\in[1,k]}\frac{1}{n-l+1}\cdot \binom{2n-2l}{n-l}\cdot 2\binom{2l-2}{l-1}=\frac{k}{n}\binom{2k}{k}\binom{2n-2k}{n-k}
\]
holds.
\end{lemma}
\textbf{Proof.} We show the identity by induction with respect to $k$.

(i) The case that $k=1$ is trivial.

(ii) Suppose that the lemma holds for $k = m < n-1$. Then we have
\begin{eqnarray*}
& &
\sum_{l\in[1,m+1]}\frac{1}{n-l+1}\cdot \binom{2n-2l}{n-l}\cdot 2\binom{2l-2}{l-1}\\
&=&
\sum_{l\in[1,m]}
\frac{1}{n-l+1}\cdot \binom{2n-2l}{n-l}\cdot 2\binom{2l-2}{l-1} 
+ \frac{1}{n - m}\cdot \binom{2n - 2m - 1}{n - m -1} \cdot 2\binom{2m}{m}.
\end{eqnarray*}
By the induction hypothesis, the above quantity is
\begin{eqnarray*}
& &
\frac{m}{n}\binom{2m}{m}\binom{2n-2m}{n-m}+\frac{1}{n-m}\cdot \binom{2n-2m-2}{n-m-1}\cdot 2\binom{2m}{m}
\\
&=&
\frac{m}{n}\binom{2m}{m}\frac{(2n-2m)(2n-2m-1)}{(n-m)^2}\binom{2n-2m-2}{n-m-1}+\frac{1}{n-m}\cdot \binom{2n-2m-2}{n-m-1}\cdot 2\binom{2m}{m}
\\
&=&
\frac{2}{n-m}\binom{2m}{m}\binom{2n-2m-2}{n-m-1} 
\left\{ \frac{m}{n}(2n-m-1)+1 \right\}
\\
&=&
\frac{2}{n-m}\binom{2m+2}{m+1}\frac{(m+1)^2}{(2m+2)(2m+1)}\binom{2n-2m-2}{n-m-1}\frac{(n-m)(2m+1)}{n}
\\
&=&
\frac{m+1}{n}\binom{2m+2}{m+1}\binom{2n-2m-2}{n-m-1}\\
&=&
\frac{m+1}{n}\binom{2(m+1)}{m+1}\binom{2n-2(m+1)}{n-(m+1)}.
\end{eqnarray*}
This is nothing but the identity for $k=m+1$. 

By (i) and (ii) the identity is proved for any $k \in [1, n-1]$. 
$\Box$

The following is the main theorem of the present paper:

\begin{theorem}\label{main theorem}
\[
\#\mathcal{P}_{2n,2k}^{+}=\frac{k}{n}\binom{2k}{k}\binom{2n-2k}{n-k}.
\]
\end{theorem}
\textbf{Proof.}
For $k=0$, it is trivial. For $k=n$, it follows from Proposition \ref{positive}.\\
The equation for the case of $0<n<k$ directly follows from Lemma \ref{path counting} and Lemma \ref{identity}.
$\Box$

As a corollary, we obtain an explanation for WHY:

\begin{theorem}\label{why}
For the simple random walk on $\mathbb{Z}$, the identity below holds.
\[
\mathbf{P}(\:X(2n) {\rm \: is\: on\: the\: positive\: side\:} |\:\: T_{2n}=2k)=\frac{k}{n}.
\]
\end{theorem}
\textbf{Proof.}
By the definition of the conditional probability and from Theorem \ref{arcsine} and Theorem \ref{main theorem}, we have
\[
\mathbf{P}(\:X(2n) \text{ is on the positive side } \ |\ T_{2n}=2k)
=\frac{\#\mathcal{P}_{2n,2k}^{+}}{\#\mathcal{P}_{2n,2k}}
=\frac{\frac{k}{n}\binom{2k}{k}\binom{2n-2k}{n-k}}{\binom{2k}{k}\binom{2n-2k}{n-k}}
=\frac{k}{n}.
\]$\Box$

This is why 
\begin{center}
``the probability that we can now find on the positive side \\
a simple random walker 
who has already spent 70 percent of time on the positive side is 70 persent.''
\end{center}
 

\section{Marchenko-Pastur Law}
For the Brownian motion $x(\tau)$, it is well known that 
\[
\mathbf{P}(\:\mu(\tau \in [0,t] \ |\ x(\tau)\geq 0)\leq rt)=\frac{1}{\pi}\int_0^{r}\frac{dx}{\sqrt{x(1-x)}},
\]
where $\mu(dx)=dx$ denotes the Lebesgue measure of the real line. 
The fact above is called the Arcsine law \cite{FEL}.

Based on the theorems in the previous section and on the standard transition from digital to analog, 
we easily obtain the arcsine law for the Brownian motion $x(\tau)$ with positive ends: 

\begin{theorem}\label{March}
We have
\[
\mathbf{P}(\:\mu(\tau \in [0,t]\:\:\ |\  \:\:x(\tau)\geq 0)\leq rt)=\frac{2}{\pi}\int_0^{r}\frac{xdx}{\sqrt{x(1-x)}}
=\frac{2}{\pi}\int_0^{r}\sqrt{\frac{x}{1-x}}dx,
\]
and dually,
\[
\mathbf{P}(\:\mu(\tau \in [0,t]\:\:\ |\  \:\:x(\tau)\leq 0)\leq rt)=\frac{2}{\pi}\int_0^{r}\frac{(1-x)dx}{\sqrt{x(1-x)}}
=\frac{2}{\pi}\int_0^{r}\sqrt{\frac{1-x}{x}}dx,
\]
where $\mu(dx)=dx$ denotes the Lebesgue measure of the real line. 
\end{theorem}
  
The probability law above is nothing but the Marchenko-Pastur law, 
which plays fundamental and  universal roles 
in the theory of random matrices and free probability. 
It strongly suggests that there is some hidden relationship between the classical probability, combinatorics, random matrices and quantum probability.

\end{document}